\documentclass[12pt]{iopart}

\newtheorem{theorem}{Theorem}
\newtheorem{corollary}{Corollary}
\newtheorem{lemma}{Lemma}
\newtheorem{remark}{Remark}
\usepackage{graphicx}

\begin{document}

\title[Cram\'{e}r moderate deviations]{Cram\'{e}r moderate deviations for the elephant random walk}

\author{X. Fan$^1$, H. Hu$^2$  and X. Ma$^1$}

\address{$^1$ Center for Applied Mathematics,
Tianjin University, Tianjin 300072, China}
\ead{fanxiequan@hotmail.com}
\address{$^2$ School of Mathematics and Statistics, Northeastern University at Qinhuangdao, Qinhuangdao, China}

\begin{abstract}
 We establish some limit theorems for the elephant random walk, including
 Berry-Esseen's bounds, Cram\'{e}r moderate deviations and local limit theorems.  These limit theorems
can be regarded as refinements of the central limit theorem  for the elephant random walk.
Moreover, by these limit theorems, we conclude that the convergence rate of normal approximations and the domain of attraction of normal distribution
mainly depend  on a memory parameter $p$ which lies between $0$ and $3/4.$  \\

\end{abstract}
%

%
%
%
%

\section{Introduction}
 The elephant random walk (ERW) was introduced by Sch{\"u}tz and Trimper \cite{schutz2004elephants}  in order to study the memory effects in the non-Markovian random walk. The model has a link to a famous saying that elephants can remember where they have been.
 Since  the seminal work of Sch{\"u}tz and Trimper \cite{schutz2004elephants}, the ERW has recently attracted a lot of attentions.
A wide range  of literature is available for the asymptotic behavior of the ERW  and its extensions, see  \cite{baur2016elephant}-\cite{CGS17} and \cite{G19}. 
Baur and Bertoin \cite{baur2016elephant} derived the functional limit theorem   via a method of connection to P\'{o}lya-type urns.
Coletti, Gava and Sch{\"u}tz \cite{coletti2017central,CGS17} proved the central limit theorem  (CLT) and a  strong invariance principle for $p \in [0,  3/4]$ and a law of large numbers for $p \in [0, 1).$ They also showed that if $p \in (3/4, 1]$, then the ERW converges to a non-degenerate random variable which is not normal.
V\'{a}zquez Guevara \cite{G19} gave the almost sure central limit theorems.
For the multi-dimensional ERW, we refer to  \cite{BL19}, where
 Bercu and Lucile have established  the CLT.

 The one-dimensional ERW can be defined as  follows.
It starts at time $n=0$, with  position $S_0=0$. At time $n=1$, the elephant
 moves to $1$ with probability $q$ and to $-1$ with probability $1-q$, where $q \in [0, 1].$
 So the position of the elephant at time $n=1$ is given by $S_1=X_1,$ with $X_1$
 a Rademacher $\mathcal{R}(q)$ random variable.
 At time $n+1,$ for $n\geq 1,$ an integer $n'$ is chosen from the set $\{1, 2,\ldots, n\}$
 uniformly at random. Then $X_{n+1}$ is determined stochastically by the following rule.
 If $X_{n'}=1$, then
 \begin{displaymath}
 X_{n+1} =   \left\{ \begin{array}{ll}
1 & \textrm{with  probability  $p$}\\
-1 & \textrm{with  probability  $1-p$}.
\end{array} \right.
\end{displaymath}
 If $X_{n'}=-1$, then
 \begin{displaymath}
 X_{n+1} =   \left\{ \begin{array}{ll}
1 & \textrm{with  probability  $1-p$}\\
-1 & \textrm{with  probability  $p$}.
\end{array} \right.
\end{displaymath}
Thus, for $n\geq1,$  the position of the elephant at time $n+1$ is
$$S_{n+1}=\sum_{i=1}^{n+1}X_i,$$
 where
$$ X_{n+1}=\alpha_nX_{\beta_n},$$
with $\alpha_n$ has a Rademacher distribution $\mathcal{R}(p)$, $p\in [0, 1],$ and $\beta_n$ is uniformly distributed over the integers $\{1, 2,\ldots, n\}$. Moreover, $\alpha_n$ is independent of $X_{1},...,X_{n}$.
Here  $p $ is called  the memory parameter. The ERW is respectively  called diffusive,  critical and   superdiffusive according to  $p \in [0, 3/4),$ $p=3/4$ and $p \in (3/4, 1]$.

Recently,   Bercu \cite{bercu2017martingale} presented the following CLT for the ERW: if  $  p \in [0,  3/4)
$, then
	\begin{equation}\label{ap01}
	\frac{S_n}{\sqrt{n /(3-4p)}}\stackrel{\mathbf{D}}{\longrightarrow} \mathcal{N}(0, 1),\ \ \ n\rightarrow \infty;
	\end{equation}
and if $p=3/4,$ then
\begin{equation}\label{ap02}
\frac{S_n}{\sqrt{n\log n}}\stackrel{\mathbf{D}}{\longrightarrow}\mathcal{N}(0,1),   \ \ \ n\rightarrow \infty,
\end{equation}
where $\stackrel{\mathbf{D}}{\longrightarrow}$ stands for convergence in distribution. He also showed that when $p \in (3/4, 1]$, the ERW converges to a non-normal  random variable.

In this paper we are interested in the absolute and relative errors of the normal approximations (\ref{ap01}) and (\ref{ap02}),
and therefore we focus on the case where the memory parameter $p$   lies between $0$ and $3/4.$
For the absolute errors of  normal approximations, we establish   Berry-Esseen's bounds, which show  that convergence rate of the absolute errors
depends on the memory parameter $p$.
For the relative errors of   normal approximations, we obtain the Cram\'{e}r  moderate deviations,
which conclude that the domain of attraction of normal distribution
mainly depends on the memory parameter $p$.  Moreover, the local limit theorems for the ERW are also established.
Notice that when $p=0$ or $1/2,$ the ERW   reduces  to the classical symmetric random walk. Thus we do not consider the  cases   $p=0$ and $1/2$.

The paper is organized as follows. In Section \ref{secm}, we present our main results,  including
   Berry-Esseen's bounds, Cram\'{e}r moderate deviations and  the local limit theorems. Section \ref{simu} contains  a simulation study and applications
   of our results.   The proofs of our results are given in Section \ref{secp}. At the last section, we give a conclusion for the paper.

Throughout the paper,  $C$ and $C_p,$ probably supplied with some indices,
denote respectively a generic positive absolute constant and a generic positive constant depending only on $p.$  For two sequences of positive numbers $(a_n)$ and $(b_n)$,   write $a_n \asymp b_n$ if there exists an absolute constant $C>0$ such that ${a_n}/{C}\leq b_n\leq C a_n$ for all sufficiently large $n$. We also  write $a_n \sim b_n$ if  $\lim_{n\rightarrow \infty} a_n/b_n =1$. Moreover, $\theta $  stands for values satisfying $\left| \theta  \right| \leq 1$.

\section{Main results}\label{secm}
\subsection{Berry-Esseen's bounds}
Denote by $\Phi(t)$  the standard normal distribution function, and
$$D\left(X \right)=\sup _{t\in\mathbf{R}} \Big|\mathbf{P}(X  \leq t)-\Phi(t) \Big|  $$
the absolute error  of the normal approximation for $X$. The upper bounds of $D(X)$ are called as   Berry-Esseen's  bounds.
For $p \in (0, 3/4],$ denote
$$ a_n=\frac{\Gamma(n)\Gamma(2p)}{\Gamma(n+2p-1)} \ \ \ \  \ \ \ \   \textrm{and} \ \ \ \ \ \ \ \    v_n=\sum_{i=1}^n a_i^2.$$
By Stirling's formula $$\log  \Gamma(x)  =(x-\frac12) \log x -x + \frac12 \log 2 \pi + O(\frac{1}{x})  \ \ \ \  \ \textrm{as}\ \ x \rightarrow \infty,$$ we deduce that
\begin{eqnarray}\label{a10}
\lim\limits_{n\to \infty} a_n n^{ 2p-1} =\Gamma(2p).
\end{eqnarray}
Moreover, in the diffusive regime $( p\in(0, 3/4))$, we have
\begin{eqnarray}\label{a15}
\lim\limits_{n\to \infty} \frac{ v_n}{ n^{3-4p}} = \frac{\Gamma{(2p)}^2}{3-4p},
\end{eqnarray}
and, in the critical regime $(p=3/4)$,  it holds
\begin{eqnarray}\label{a16}
\lim\limits_{n\to \infty} \frac{ v_n}{  \log n}    =\frac{\pi}{4}.
\end{eqnarray}

Our first result concerns with  Berry-Esseen's  bounds for the  ERW, which shows that the absolute errors of normal approximations  mainly depend on the memory parameter $p$.
\begin{theorem}\label{thm1}
 Assume $  p  \in (0,   3/4]$ and $p\neq 1/2$.
The following Berry-Esseen's bounds hold.
\begin{description}
  \item[\textbf{[i]}]  If $  p \in (0, 1/2),$  then
	\begin{equation}\label{ineq2}
	D\Big(  \frac{a_n S_n}{ \sqrt{v_n}}   \Big)\leq    C  \frac{\log n }{ \sqrt{n}  }  .
	\end{equation}
 \item[\textbf{[ii]}]  If $   p \in (1/2, 3/4),$  then
	\begin{equation}\label{ineq3}
	D\Big(  \frac{a_n S_n}{ \sqrt{v_n}}   \Big)\leq    C_p  \frac{\log n }{    n^{ (3-4p) /2 } \ } .
	\end{equation}
  \item[\textbf{[iii]}] If $p=3/4,$ then
   \begin{equation}\label{ineq4}
   D\Big(  \frac{a_n S_n}{ \sqrt{v_n}}   \Big)\leq C \frac{ \log \log n     }{\sqrt{\log n}  }.
   \end{equation}	
\end{description}
\end{theorem}	


%

The following corollary is a simple consequence of Theorem \ref{thm1}, which bridges   Theorem \ref{thm1} and the CLT of  Bercu \cite{bercu2017martingale}.

\begin{corollary}\label{cor01}
 Assume $  p  \in (0,   3/4]$ and $p\neq 1/2$.
The following Berry-Esseen's bounds hold.
\begin{description}
  \item[\textbf{[i]}]  If $ p \in (0, 1/2),$  then
	\begin{equation}\label{ineqsf11}
	D\Big(  \frac{ S_n}{ \sqrt{n /(3-4p)}}   \Big)\leq   C     \bigg(\frac{\log   n }{ \sqrt{n}  }   + \Big|\frac{\sqrt{v_n}}{a_n \sqrt{n/(3-4p)}}  -1 \Big| \bigg)  .
	\end{equation}
     \item[\textbf{[ii]}]  If $  p \in (1/2, 3/4),$  then
	\begin{equation}\label{ineqsf12}
	D\Big(  \frac{ S_n}{ \sqrt{n /(3-4p)}}   \Big)\leq   C_p     \bigg(\frac{\log n }{   n^{(3-4p)/2 } \ }  + \Big|\frac{\sqrt{v_n}}{a_n \sqrt{n/(3-4p)}}  -1 \Big| \bigg)  .
	\end{equation}
  \item[\textbf{[iii]}] If $p=3/4,$ then
   \begin{equation}\label{ineqsf2}
   D\Big(  \frac{  S_n}{ \sqrt{ n \log n}}   \Big)\leq  C \bigg(  \frac{ \log \log n     }{\sqrt{\log n}  } +   \Big|\frac{\sqrt{v_n}}{a_n\sqrt{n\log n}}  -1 \Big|\bigg).
   \end{equation}	
\end{description}
\end{corollary}

 By the equalities  (\ref{a10}) and (\ref{a15}), it is easy to verify  that   for $0< p < 3/4,$  it holds
  \begin{equation}\label{idf2}
  \displaystyle \lim_{n\rightarrow \infty} \frac{\sqrt{v_n} } { a_n\sqrt{n/(3-4p)}}  =1.
  \end{equation}
 Similarly,  for $ p = 3/4,$ by the equalities  (\ref{a10}) and (\ref{a16}) and the fact $\Gamma(3/2)=\sqrt{ \pi}/2$, it holds  $\displaystyle \lim_{n\rightarrow \infty} \frac{\sqrt{v_n} }{ a_n\sqrt{n \log n}}   =1.$   Thus the results  of Corollary \ref{cor01} coincide  with the CLT of  Bercu \cite{bercu2017martingale}, that is, (\ref{ap01}) and (\ref{ap02}).

\subsection{Cram\'{e}r   moderate deviations}
	The following theorem gives some Cram\'{e}r  moderate deviations for the ERW.
\begin{theorem}\label{thmg2}
 Assume $  p  \in (0,   3/4]$ and $p\neq 1/2$.
 We have the following equalities.
 \begin{description}
  \item[\textbf{[i]}]  If $p \in (0,  1/2)$,  then there is an absolute constant $\alpha_0 >0$ such that for all $0\leq x \leq \alpha_0 \, \sqrt{n},$
	\begin{equation}\label{fth1}
\frac{\mathbf{P}(a_n S_n/ \sqrt{v_n} \geq x)}{1-\Phi \left( x\right)}=\exp\bigg\{\theta C_{\alpha_0 } \bigg(\frac{x^3}{\sqrt{n}}   +  (1+x) \frac{ \log n }{\sqrt{n}}  \bigg)  \bigg\}.
\end{equation}
  \item[\textbf{[ii]}]  If $p \in (1/2, 3/4)$,  then there is an absolute constant $\alpha_0 >0$ such that for all $0\leq x \leq \alpha_0 \,     n^{ (3-4p) /2 },$
	\begin{equation}
\frac{\mathbf{P}(a_n S_n/ \sqrt{v_n} \geq x)}{1-\Phi \left( x\right)}=\exp\bigg\{\theta C_{\alpha_0, p} \bigg(\frac{x^3}{    n^{ (3-4p) /2 } }  +  (1+x) \frac{ \log n }{   n^{ (3-4p) /2 } } \bigg)  \bigg\}.
\end{equation}

  \item[\textbf{[iii]}] If $ p=3/4$,  then there is an absolute constant $\alpha_0 >0$ such that for all $0\leq x \leq \alpha_0 \, \sqrt{\log n},$
	\begin{equation}
\!\!\! \frac{\mathbf{P}(a_n S_n/ \sqrt{v_n} \geq x)}{1-\Phi \left( x\right)}=\exp\bigg\{\theta C_{\alpha_0 } \bigg(\frac{x^3}{\sqrt{ \log n}}  +  (1+x) \frac{ \log \log n  }{\sqrt{\log n}}  \bigg)  \bigg\}.
\end{equation}
\end{description}
Moreover, the same equalities hold when $\frac{\mathbf{P}(a_n S_n/ \sqrt{v_n} \geq x)}{1-\Phi \left( x\right)}$ is replaced by $\frac{\mathbf{P}(a_n S_n/ \sqrt{v_n} \leq -x)}{\Phi \left( -x\right) }$.
\end{theorem}

It is easy to see that $|e^x-1| \leq e |x|$ for all $|x|\leq 1.$
By Theorem \ref{thmg2}, we have the following corollary, which shows that  the domain of attraction of normal distribution
mainly depends on the memory parameter $p$.

\begin{corollary}\label{co02}
 Assume $  p  \in (0,   3/4]$ and $p\neq 1/2$. The following results   hold.
  \begin{description}
  \item[\textbf{[i]}]  If $  p \in (0, 1/2)$, then for all $0\leq x \leq   n^{1/6} ,$
  \begin{equation}
 \bigg| \frac{\mathbf{P}(a_n S_n/ \sqrt{v_n} \geq x)}{1-\Phi \left( x\right)} - 1 \bigg| \leq C  \bigg(\frac{x^3}{\sqrt{n}}   +  (1+x) \frac{ \log n }{\sqrt{n}}  \bigg).
\end{equation}
In particular, the last  inequality implies that
\begin{equation}\label{sgfdf23}
\frac{\mathbf{P}(a_n S_n/ \sqrt{v_n} \geq x)}{1-\Phi \left( x\right)}=1+o(1)
\end{equation}
uniformly for  $0 \leq x =o(n^{1/6})$ as $n\rightarrow \infty$.
  \item[\textbf{[ii]}]  If $ p \in (1/2, 3/4)$, then for all $0\leq x \leq   n^{(3-4p)/6} ,$
  \begin{equation}
 \bigg| \frac{\mathbf{P}(a_n S_n/ \sqrt{v_n} \geq x)}{1-\Phi \left( x\right)} - 1 \bigg| \leq C_p  \bigg(\frac{x^3}{    n^{ (3-4p) /2 } }  +  (1+x) \frac{ \log n }{   n^{ (3-4p) /2 }}  \bigg).
\end{equation}
In particular, the last  inequality implies that (\ref{sgfdf23}) holds
uniformly for  $0 \leq x =o(  n^{(3-4p)/6}) $ as $n\rightarrow \infty$.

  \item[\textbf{[iii]}]  If $ p=3/4$, then for all $0\leq x \leq  (\log n)^{1/6} ,$
  \begin{equation}
 \bigg| \frac{\mathbf{P}(a_n S_n/ \sqrt{v_n} \geq x)}{1-\Phi \left( x\right)} - 1 \bigg| \leq C \bigg(\frac{x^3}{\sqrt{ \log n}}  +  (1+x) \frac{   \log \log n}{\sqrt{\log n}}  \bigg).
\end{equation}
In particular, the last   inequality implies that (\ref{sgfdf23}) holds
uniformly for  $0 \leq x =o( (\log n)^{1/6}) $ as $n\rightarrow \infty$.
\end{description}
Moreover, the same equalities hold when $\frac{\mathbf{P}(a_n S_n/ \sqrt{v_n} \geq x)}{1-\Phi \left( x\right)}$ is replaced by $\frac{\mathbf{P}(a_n S_n/ \sqrt{v_n} \leq -x)}{\Phi \left( -x\right) }$.
\end{corollary}

From Theorem \ref{thmg2}, we can deduce the following moderate deviation principles (MDP) for the ERW.
\begin{corollary}\label{co03}
 Assume $  p  \in (0,   3/4]$ and $p\neq 1/2$. The following MDP   hold.
 \begin{description}
  \item[\textbf{[i]}] Assume   $p \in (0, 1/2)$.
Let $b_n$ be any sequence of real numbers satisfying $b_n \rightarrow \infty$ and $b_n/\sqrt{n}\rightarrow 0$
as $n\rightarrow \infty$.  Then  for each Borel set $B$,
\begin{eqnarray}
- \inf_{x \in B^o}\frac{x^2}{2} &\leq & \liminf_{n\rightarrow \infty}\frac{1}{b_n^2}\ln \mathbf{P}\bigg(\frac{   a_n S_n }{b_n \sqrt{v_n} }  \in B \bigg) \nonumber \\
 &\leq& \limsup_{n\rightarrow \infty}\frac{1}{b_n^2}\ln \mathbf{P}\bigg(\frac{   a_n S_n }{b_n \sqrt{v_n} }   \in B \bigg) \leq  - \inf_{x \in \overline{B}}\frac{x^2}{2} \, ,   \label{MDP}
\end{eqnarray}
where $B^o$ and $\overline{B}$ denote the interior and the closure of $B$, respectively.
  \item[\textbf{[ii]}]  Assume   $ p \in (1/2, 3/4)$.
Let $b_n$ be any sequence of real numbers satisfying $b_n \rightarrow \infty$ and $b_n/ n^{(3-4p)/2}  \rightarrow 0$
as $n\rightarrow \infty$.  Then (\ref{MDP}) holds.

  \item[\textbf{[iii]}]  Assume   $  p = 3/4$.
Let $b_n$ be any sequence of real numbers satisfying $b_n \rightarrow \infty$ and $b_n/\sqrt{\ln n} \rightarrow 0$
as $n\rightarrow \infty$.  Then (\ref{MDP}) holds.
\end{description}

\end{corollary}

\subsection{Local limit theorems}
The following theorem gives  asymptotic expansions for the probability   $\mathbf{P}( S_n=k  )$.
\begin{theorem} \label{lclt}
  Assume $  p  \in (0,   3/4]$ and $p\neq 1/2$. The following results hold.
 \begin{description}
  \item[\textbf{[i]}]  If $  p\in (0, 1/2)$,  then   for all $ |k|= o(n^{2/3}),$
	\begin{equation} \label{f20}
\Bigg|\frac{\mathbf{P}( S_n=k  )}{ \frac{a_n}{\sqrt{2 \pi v_n}   } \exp\{ - \frac{(a_nk)^2}{2 v_n}  \}  } - 1 \Bigg| \leq  C \bigg(\frac{|k|}{n}+ \frac{k^2}{ n^{3/2}} +\frac{\log n}{\sqrt{n}} \bigg).
\end{equation}
In particular, the last inequality implies   the following local limit theorem:
\begin{equation}\label{locclt}
 \frac{\mathbf{P}( S_n=k  )}{ \frac{a_n}{\sqrt{2 \pi v_n}   } \exp\{ - \frac{(a_nk)^2}{2 v_n} \}  }    = 1+ o(1)
\end{equation}
uniformly for   $ |k|= o(n^{2/3})$ as $n\rightarrow \infty$.
  \item[\textbf{[ii]}]  If $   p \in(1/2,   3/4)$,  then   for all $ |k|= o(n^{(3-2p)/3}),$
	\begin{equation} \label{f20}
\Bigg|\frac{\mathbf{P}( S_n=k  )}{ \frac{a_n}{\sqrt{2 \pi v_n}   } \exp\{ - \frac{(a_nk)^2}{2 v_n}  \}  } - 1 \Bigg| \leq  C_p \bigg(\frac{|k|}{n}+ \frac{k^2}{ n^{(5-4p)/2 }} +\frac{\log n}{   n^{ (3-4p) /2 }} \bigg).
\end{equation}
In particular, the last inequality implies that  (\ref{locclt}) holds
uniformly for   $ |k|= o(n^{(3-2p)/3})$ as $n\rightarrow \infty$.
  \item[\textbf{[ii]}] If $ p=3/4$,  then   for all $ |k|= o((n\log n)^{1/2} ),$
	\begin{equation}\label{f21}
\Bigg|\frac{\mathbf{P}( S_n=k  )}{ \frac{a_n}{\sqrt{2 \pi v_n}   } \exp\{ - \frac{(a_nk)^2}{2 v_n} \}  } - 1 \Bigg| \leq    C   \bigg(\frac{|k|}{n \log n}+ \frac{k^2}{ n (\log n)^{3/2} } +  \frac{\log \log n}{\sqrt{\log n}}\bigg).
\end{equation}
Thus  (\ref{locclt}) holds uniformly for   $ |k|= o((n\log n)^{1/2} )$ as $n\rightarrow \infty$.
\end{description}
\end{theorem}

Denote by
$$L\left(S_n \right)=  \sup_{k \in \mathbf{Z}  } \bigg  |\mathbf{P}(S_n =k  )- \frac{a_n}{\sqrt{2 \pi v_n}}\exp\Big\{-  \frac{(a_n k)^2}{2 v_n}  \Big\} \bigg| $$
the absolute error  for  the local limit theorem.
\begin{corollary}\label{th025}
 Assume $  p  \in (0,   3/4)$ and $p\neq 1/2$. The following inequalities hold for $L\left(S_n \right)$.
 \begin{description}
  \item[\textbf{[i]}]  If $  p \in (0, 1/2)$,  then
	\begin{equation}
L\left(S_n \right) \leq  C \frac{\log n}{n}  .
\end{equation}
\item[\textbf{[ii]}]    If $ p \in (1/2, 3/4)$,  then
	\begin{equation}
L\left(S_n \right) \leq  C_p \frac{\log n}{ n^{2-2p} }  .
\end{equation}
\end{description}
\end{corollary}

From the last corollary, it is easy to see that the absolute errors for the local limit theorem
are much smaller than that for the CLT (cf.\ Berry-Esseen's bounds in Theorem \ref{thm1}).


\section{A simulation study and applications}\label{simu}
\subsection{A simulation study}
\begin{figure}
\includegraphics[width=0.9\linewidth]{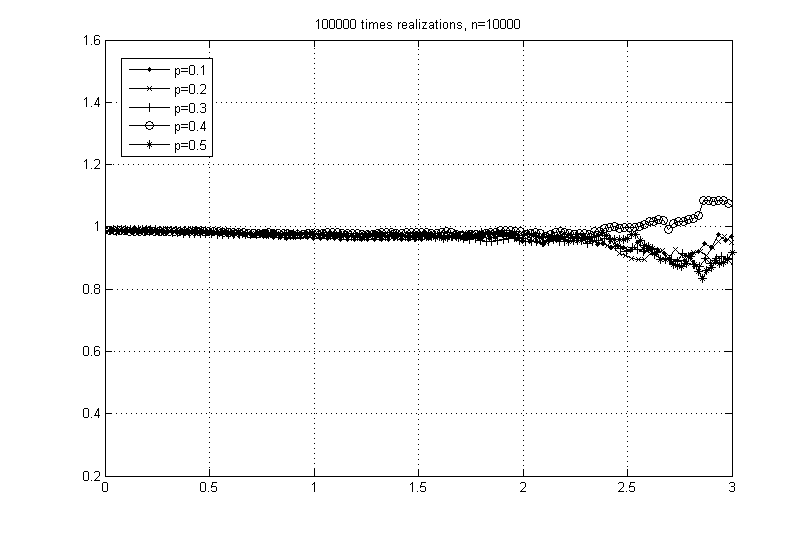}
\caption[]{Ratios  $\frac{\mathbf{P}(a_n S_n/ \sqrt{v_n} \geq x)}{1-\Phi \left( x\right)}$ for ERW
with $n=10^4$    and $p = 0.1, 0.2,..., 0.5.$    }
\label{ssf1}
\vspace{0.5cm}

\includegraphics[width=0.9\linewidth]{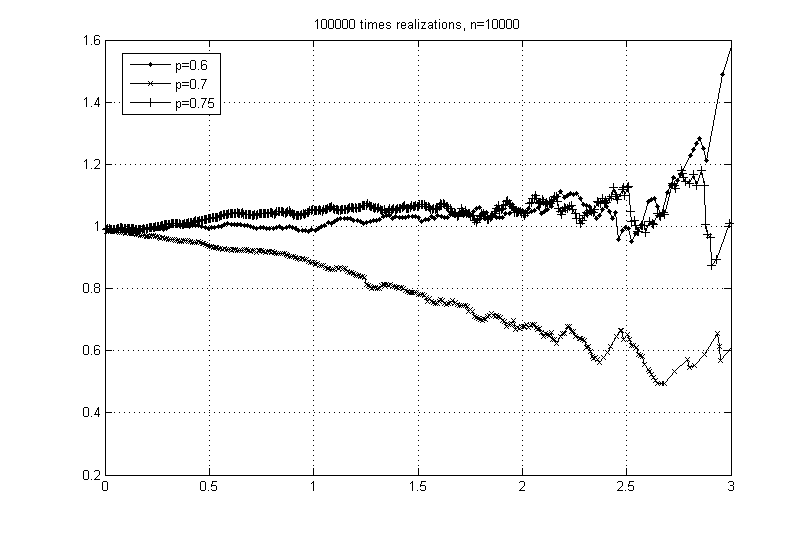}
\caption[]{ Ratios  $\frac{\mathbf{P}(a_n S_n/ \sqrt{v_n} \geq x)}{1-\Phi \left( x\right)}$ for ERW
with $n=10^4$ and $p = 0.6$, $0.7.$  \\
For $p=0.75,$ we  adopt the ratio $\frac{\mathbf{P}(S_n/  \sqrt{n\log n}   \geq  x)}{1-\Phi \left( x\right)}$.  }
\label{ssf2}
\end{figure}

In this section, we study the normal approximation accuracies for the ERW.
We let $n=10^4$, and choose 8 levels of $p$:\ $p = 0.1, 0.2,..., 0.7, 0.75.$
Figures \ref{ssf1} and   \ref{ssf2} show  the Cram\'{e}r moderate deviation ratios  $ \frac{\mathbf{P}(a_n S_n/ \sqrt{v_n} \geq x)}{1-\Phi(x)} $ and $\frac{\mathbf{P}(S_n/  \sqrt{n\log n}   \geq  x)}{1-\Phi \left( x\right)}$,
where the probabilities  $\mathbf{P}(a_n S_n/ \sqrt{v_n} \geq x)$ and $\mathbf{P}(S_n/  \sqrt{n\log n}   \geq  x)$ are approximated by simulating
$10^5$ realizations of the ERW.
When $p=0.1, 0.2,...,   0.5$,  we see that the Cram\'{e}r moderate deviation ratios have comparable performance; see Figure \ref{ssf1}.
When $p=0.6, 0.7, 0.75,$ the normal approximations are worse than the case of $p=0.1, 0.2,...,0.5;$ see Figure \ref{ssf2}.
Moreover, as $x$ moves away from $0$, the moderate deviation approximations
  become worse. These simulation results coincide with  Corollary \ref{co02}.

 \subsection{Applications}
  Berry-Esseen's  bounds and Cram\'{e}r moderate deviations for the ERW can be applied to construct confidence  lower limit for
the   memory parameter $p$. Assume   $p \in (0, 3/4)$ and $\kappa  \in (0,  1)$.
By the inequalities (\ref{ineqsf11}) and (\ref{ineqsf12}), then we have
 \begin{equation}\label{iddsfs02}
 \mathbf{P}\bigg( | S_n/  \sqrt{n /(3-4p)}  |  \leq \Phi^{-1}(1-\kappa /2) \bigg) \rightarrow  1- \kappa , \ \ \ \ n \rightarrow \infty.
 \end{equation}
Clearly, $| S_n/  \sqrt{n /(3-4p)}  |  \leq \Phi^{-1}(1-\kappa /2)$  means that $p\geq \frac{1}{4}\big(3-n\big(  \frac{\Phi^{-1}(1-\kappa /2)}{S_n} \big)^2  \big).$
Thus (\ref{iddsfs02}) implies that
$
  \frac{1}{4}\big(3-n\big(  \frac{\Phi^{-1}(1-\kappa /2)}{S_n} \big)^2  \big)
$
is a  $1-\kappa $ lower confidence limit for $p$, for $n$ large enough. When $\kappa$ is replaced by $\kappa_n$ such that $\kappa_n \rightarrow 0 $ and $| \log \kappa_n|=o(\sqrt{ n})$ as $n\rightarrow \infty$, similar results are allowed to be established via  Cram\'{e}r moderate deviations
for the ERW. Conversely,  if $p$ is known, Cram\'{e}r moderate deviations can also be used to
interval estimations for the position of the elephant in the ERW model.

\section{Proofs of Theorems} \label{secp}
\subsection{Preliminary lemmas}
Set $$\gamma _n=1+\frac{2p-1}{n} $$ and put $a_1=1.$  For $n\geq 2,$ it is easy to see that
$$a_n=\frac{\Gamma(n)\Gamma(2p)}{\Gamma(n+2p-1)}=\prod_{i=1}^{n-1}\frac{1}{\gamma_i}.$$
Define the filtration $\mathcal{F}_n=\sigma\{X_i: 1\leq i\leq n\}  $ and
\begin{equation}\label{mndef}
M_n=a_nS_n \qquad a.s.
\end{equation}
It is easy to verify that  $(M_n, \mathcal{F}_n)_{n\geq1}$ is a martingale. Indeed,  for
all $n\geq 1,$ we have
\begin{eqnarray}
\mathbf{E}[M_{n+1}| \mathcal{F}_n]&=&\mathbf{E}[a_{n+1}(S_n+\alpha_nX_{\beta_n})|\mathcal{F}_n]\nonumber\\
&=&a_{n+1}S_n+a_{n+1}\mathbf{E}[\alpha_n]\mathbf{E}[X_{\beta_n}|\mathcal{F}_n]\nonumber\\
&=&a_{n+1} \Big(S_n  + (2p-1) \frac{S_n}{n} \Big) \nonumber\\
&=&a_nS_n\nonumber\\
&=&M_n \qquad a.s.
\end{eqnarray}
Moreover, we can rewrite $(M_n)_{n\geq1}$ in the following additive form
\begin{equation}\label{fsdfsdf}
M_n=\sum_{i=1}^{n}a_i\varepsilon_i,
\end{equation}
where $\varepsilon_i=S_i-\gamma_{i-1}S_{i-1}$ with $S_0=0.$
 Let $(\Delta M_n)_{n\geq1}$ be the martingale differences defined by $\Delta M_1=M_1$ and  for $n\geq 2,$
$$\Delta M_n=M_n-M_{n-1}.$$
Denote by $\langle M\rangle_n $   the quadratic variation $$\langle M\rangle_n=\sum_{i=1}^n\mathbf{E}[\Delta M_i^2|\mathcal{F}_{i-1}].$$

In the proofs of theorems, we need the following two lemmas for the boundness of martingale differences and  the
convergence of  quadratic variation.
\begin{lemma}\label{lem1}
	For each $n\geq1$ and $p \in [0, 1]$, it holds $$|\Delta M_n|\leq 2a_n.$$
\end{lemma}
\noindent\emph{Proof.}
It holds obviously for $n=1$.
Observe that
\begin{eqnarray*}
	\Delta M_n = a_nS_n-a_{n-1}S_{n-1}
	 = a_nX_n-a_n\frac{S_{n-1}}{n-1} (2p-1).
\end{eqnarray*}
Since $|X_n| \leq 1$, we have   $|S_{n-1}|\leq n-1$.  It is easy to see that $|\Delta M_n|\leq 2a_n.$

\begin{lemma}\label{lem2}  Assume $ p \in (0,   3/4].$ For all $1\leq i \leq n$,  we have
	 $$\bigg\|  \frac{\Delta M_i }{\sqrt{v_n}} \bigg\|_\infty \leq   \frac{2\, a_i }{\sqrt{v_n} } $$
and
\begin{eqnarray}\label{thineq31}
\displaystyle
 \bigg\|\frac{\langle M\rangle_n}{v_n}-1\bigg\|_\infty    \leq  \left\{ \begin{array}{ll}\displaystyle
\frac{ C  }{3-4p}  \frac{1}{n}  ,\ \   & \textrm{if $0< p <3/4  $,}\\
\\
\displaystyle \frac{C  }{\log n }  , \ \   & \textrm{if $  p=3/4$.}
\end{array} \right.
\end{eqnarray}
\end{lemma}
\noindent\emph{Proof.}
Clearly,   from Lemma \ref{lem1}, we have  $$\|\Delta M_i/\sqrt{v_n} \|_\infty \leq 2 a_i/\sqrt{v_n},$$
which gives the first desired inequality.
Next we give an estimation of $ \| \langle M\rangle_n/ v_n -1 \|_\infty $ for   $0<p \leq 3/4.$
From (\ref{fsdfsdf}), we get $\Delta M_k = a_k\varepsilon_k=a_k(S_k-\gamma_{k-1}S_{k-1}). $ Thus, it holds
\begin{eqnarray*}
	\mathbf{E}[(\Delta M_k)^2   |\mathcal{F}_{k-1}]&=&a_k^2 \mathbf{E}[ (S_k-\gamma_{k-1}S_{k-1})^2 |  \mathcal{F}_{k-1}]\\
&=&a_k^2 \big( \mathbf{E}[S_k^2|  \mathcal{F}_{k-1}]  -2\gamma_{k-1}S_{k-1}\mathbf{E}[ S_k |  \mathcal{F}_{k-1}]+\gamma_{k-1}^2S_{k-1}^2 \big) .
\end{eqnarray*}
It is easy to see that
\begin{eqnarray*}
  \mathbf{E}[S_k^2|  \mathcal{F}_{k-1}] &=& \mathbf{E}[ (S_{k-1}+\alpha_kX_{\beta_k})^2|\mathcal{F}_{k-1}] \\
  &=& S_{k-1}^2+ 2 S_{k-1} \mathbf{E}[  \alpha_k X_{\beta_k} |\mathcal{F}_{k-1}] +1 \\
  &=&S_{k-1}^2+ 2 \frac{2p-1}{k-1}  S_{k-1}^2   +1 \\
  &=&  (2 \gamma_{k-1} -1) S_{k-1}^2   +1
\end{eqnarray*}
and
\begin{eqnarray*}
\mathbf{E}[ S_k |  \mathcal{F}_{k-1}] &=& \mathbf{E}[  S_{k-1}+\alpha_kX_{\beta_k} |\mathcal{F}_{k-1}]
   =  S_{k-1} +   \frac{2p-1}{k-1}  S_{k-1} =  \gamma_{k-1}   S_{k-1}.
\end{eqnarray*}
Thus, we have $\mathbf{E}[(\Delta M_1)^2   ] =1=a_1^2 $ and for $k\geq2,$
\begin{eqnarray*}
	\mathbf{E}[(\Delta M_k)^2   |\mathcal{F}_{k-1}]
&=&a_k^2 \big( (2 \gamma_{k-1} -1) S_{k-1}^2   +1   -2\gamma_{k-1}^2S_{k-1}^2 +\gamma_{k-1}^2S_{k-1}^2 \big) \\
&=&a_k^2 \big( 1-( \gamma_{k-1} -1)^2 S_{k-1}^2 ) \\
&=& a_k^2   - ( 2p -1)^2 a_k^2 ( \frac{S_{k-1}}{k-1} )^2   .
\end{eqnarray*}
Hence, by the definition of $v_n$ and $M_k$, we obtain
\begin{eqnarray*}
	\langle M\rangle_n&=&v_n-(2p-1)^2\bigg(\sum_{k=1}^{n-1}\Big (\frac{a_{k+1}}{a_k}\Big)^2\Big(\frac{M_k}{k}\Big)^2\bigg).
\end{eqnarray*}
Since $\frac{a_{n+1}}{a_n}\sim 1$ as $n\rightarrow \infty$ (cf.\ (\ref{a10})),   we have
\begin{equation}
\|\langle M\rangle_n-v_n  \|_\infty \leq C_1(2p-1)^2\Big\|\sum_{k=1}^{n-1}(\frac{M_k}{k})^2\Big\|_\infty \leq C_2\sum_{k=1}^{n-1}\frac{1}{k^2}\|M_k\|_{\infty}^2.
\end{equation}
Using Lemma  \ref{lem1},   we derive that
\begin{eqnarray*}
	\|M_k\|_{\infty}^2\leq  \sum_{l=1}^k\|\Delta M_l\|_{\infty}^2 \leq 4 v_k.
\end{eqnarray*}
In the diffusive regime $0< p<3/4 $, by (\ref{a15}), we get
\begin{equation}
\|\langle M\rangle_n-v_n  \|_\infty\leq 4C_2\sum_{k=1}^{n-1}\frac{1}{k^2}v_k \leq C_3  \frac{\Gamma{(2p)}^2}{3-4p}\sum_{k=1}^{n-1}k^{1-4p}\leq  \frac{ C_4 }{3-4p}  n^{2-4p}.
\end{equation}
In the critical regime $p=3/4,$ by (\ref{a16}), we have
\begin{equation}
\|\langle M\rangle_n-v_n  \|_\infty\leq C_5 \sum_{k=1}^{n-1}\frac{\log k}{k^2} \leq   C_6  .
\end{equation}
Consequently, again by (\ref{a15}) and (\ref{a16}), we obtain (\ref{thineq31}).
This completes the proof of lemma.

\subsection{Proof of Theorem \ref{thm1}}
In the proof of Theorem \ref{thm1}, we will make use  of the following lemma of Fan \cite{XQF19}, which
gives  an exact  Berry-Esseen's bound  for martingales.
\begin{lemma}\label{lemma1}
Assume that   there exist positive numbers $\epsilon_n>0$ and $\delta_n\geq 0$, such that
for all $1\leq i\leq n,$
\begin{equation} \label{cond3}
 \|\Delta M_{i}/\sqrt{v_n}\|_\infty \leq \epsilon_n
\end{equation}
and
\begin{equation}
 \big\| \langle M \rangle_n /v_n-1\big\|_\infty   \leq \delta_n^2  \ \  \ \ \textrm{a.s.}
\end{equation}
If $\epsilon_n,  \delta_n \rightarrow 0$ as $n\rightarrow \infty,$ then
\begin{equation}
D\Big( \frac{M _n}{\sqrt{v_n}} \Big)\leq C  \,   \Big(  \epsilon_n |\log \epsilon_n |   + \delta_n   \Big) ,
\end{equation}
where   $C$ is a   positive absolute  constant.
\end{lemma}

Now we are in position to prove Theorem \ref{thm1}.
 Clearly, we have
 $$\frac{a_n S_n}{\sqrt{v_n}}  =  \frac{M_n  }{\sqrt{v_n}} = \sum_{i=1}^n   \frac{ \Delta M_i }{\sqrt{v_n}} ,$$
 and
 $(\Delta M_i/\sqrt{v_n} , \mathcal{F}_{i})_{i=1,...,n}$ is a finite sequence of martingale differences.
From Lemma \ref{lem2}, we have  $$\|\Delta M_i/\sqrt{v_n} \|_\infty \leq 2 \max_{1\leq i \leq n}a_i/\sqrt{v_n}=: \epsilon_n .$$
Using the inequalities (\ref{a10})-(\ref{a16}), we deduce that
\begin{eqnarray}\label{ineq37}
\displaystyle
\epsilon_n   \asymp  \left\{ \begin{array}{ll}\displaystyle
  \ n^{-1/2},\ \   & \textrm{if $0<p <1/2  $,}\\
\sqrt{ 3-4p  } \ n^{- (3-4p)/2},\ \   & \textrm{if $1/2 < p <3/4  $,} \\
\displaystyle    ( \log n)^{-1/2}, \ \   & \textrm{if $  p=3/4$.}
\end{array} \right.
\end{eqnarray}
Moreover, from Lemma \ref{lem2},  we have
\begin{eqnarray}\label{ineq38}
\displaystyle
 \bigg\|\frac{\langle M\rangle_n}{v_n}-1\bigg\|_\infty    \leq  \delta_n^2 := \left\{ \begin{array}{ll}\displaystyle
\frac{ C  }{3-4p}     n^{-1},\ \   & \textrm{if $0< p <3/4  $,}\\
C   \big( \log n    \big)^{-1}, \ \   & \textrm{if $  p=3/4$.}
\end{array} \right.
\end{eqnarray}
 Applying Lemma \ref{lemma1} to $M_n/\sqrt{v_n}$, we obtain the desired results.

\begin{remark}
Under the conditions  of Lemma \ref{lemma1}, El Machkouri and  Ouchti \cite{Mac07}   had obtained the following Berry-Esseen bound
$$D\big( \frac{M _n}{\sqrt{v_n}} \big)\leq C   \Big( \epsilon_n    \log n   + \delta_n \Big).$$
  For $0< p < 1/2,$ we have $\epsilon_n \asymp 1/\sqrt{n}, $ and thus for $0< p < 1/2,$ Theorem \ref{thm1}  can also be obtained via the Berry-Esseen  bound of  El Machkouri and  Ouchti \cite{Mac07} stated before.
\end{remark}


\subsection{Proof of Corollary \ref{cor01}}
  We only present the proof of Corollary \ref{cor01} for $0< p < 1/2$. For $1/2 <   p \leq 3/4$,  Corollary \ref{cor01}  can be proved in a similar way.  Clearly,  it holds
	\begin{eqnarray*}
&& D\Big(  \frac{ S_n}{ \sqrt{n/(3-4p)}}   \Big) \\
 &&= \sup_{t \in \mathbf{R}}   \Big|\mathbf{P}\Big( \frac{ S_n}{ \sqrt{n/(3-4p)}}  \leq t \Big)-\Phi(t) \Big| \nonumber \\
&&= \sup_{t \in \mathbf{R}}   \Big|\mathbf{P}\Big( \frac{ S_n}{ \sqrt{n/(3-4p)}}  \leq \frac{\sqrt{v_n}}{a_n \sqrt{n/(3-4p)}}t \Big)-\Phi\Big(\frac{\sqrt{v_n}}{a_n \sqrt{n/(3-4p)}}t\Big) \Big| \nonumber \\
&&\leq  \sup_{t \in \mathbf{R}}  \bigg|\mathbf{P}\Big( \frac{ S_n}{ \sqrt{n/(3-4p)}}  \leq \frac{\sqrt{v_n}}{a_n \sqrt{n/(3-4p)}}t \Big)-\Phi \big(t \big)  \bigg| \nonumber  \\
 &&+ \sup_{t \in \mathbf{R}} \bigg|\Phi \big(t \big) -\Phi \Big(\frac{\sqrt{v_n}}{a_n \sqrt{n/(3-4p)}}t \Big) \bigg| \nonumber \\
&&=  \sup_{t \in \mathbf{R}}  \bigg|\mathbf{P}\Big( \frac{ a_n S_n}{ \sqrt{v_n}  }  \leq  t \Big)-\Phi \big(t \big)  \bigg|
  + \sup_{t \in \mathbf{R}} \bigg|\Phi \big(t \big) -\Phi \Big(\frac{\sqrt{v_n}}{a_n \sqrt{n/(3-4p)}}t \Big) \bigg| .
	\end{eqnarray*}
By inequality (\ref{ineq2}), we get
\begin{eqnarray*}
D\Big(  \frac{ S_n}{ \sqrt{n/(3-4p)}}   \Big)
& \leq& C_1    \frac{\log n }{ \sqrt{n}  } + C_2 \Big|\frac{\sqrt{v_n}}{a_n \sqrt{n/(3-4p)}}  -1 \Big|\nonumber \\
&\leq&  C   \bigg( \frac{\log n }{ \sqrt{n}  } +  \Big|\frac{\sqrt{v_n}}{a_n \sqrt{n/(3-4p)}}  -1 \Big| \bigg),
	\end{eqnarray*}
which gives the desired inequality.

\subsection{Proof of Theorem \ref{thmg2}}

In the proof of Theorem \ref{thmg2}, we shall   make use  of the following lemma of Fan, Grama and Liu \cite{FGL20}, which gives
a Cram\'{e}r type moderate deviation for martingales. The lemma is a simple consequence of Theorem 1  of \cite{FGL20}, with
$\xi_i=  \Delta M_i / \sqrt{v_n} ,$ $\rho=1$ and $\varepsilon_n= e \epsilon_n.$
See also Grama and Haeusler \cite{GH00} for an earlier  result.
\begin{lemma}\label{lemma3}
Assume the conditions of Lemma \ref{lemma1}.
 Then there is an absolute constant $\alpha_0 >0$ such that
for all $0\leq x \leq \alpha_0 \,  \epsilon_n^{-1}, $
$$
  \frac{\mathbf{P}( M _n/\sqrt{v_n}  \geq x)}{1-\Phi \left( x\right)}=\exp\bigg\{\theta C_{\alpha_0} \Big( x^3 \epsilon_n  + x^2 \delta_n^2 + (1+ x) \left( \epsilon_n \left| \log  \epsilon_n
 \right|+  \delta_n \right)  \Big) \bigg\}
$$
and
$$
\frac{\mathbf{P}( M_n/\sqrt{v_n}  \leq -x)}{\Phi \left( -x\right) } = \exp\bigg\{\theta C_{\alpha_0} \Big( x^3 \epsilon_n  + x^2 \delta_n^2 + (1+ x) \left( \epsilon_n \left| \log  \epsilon_n
 \right|+  \delta_n \right)  \Big) \bigg\}.
$$
\end{lemma}

Now we are in position to prove Theorem \ref{thmg2}.
Recall that
$$\frac{a_n S_n}{\sqrt{v_n}}  =  \frac{M_n  }{\sqrt{v_n}} = \sum_{i=1}^n   \frac{ \Delta M_i }{\sqrt{v_n}}.$$
According to the proof of  Theorem \ref{thm1}, the  conditions of Lemma \ref{lemma1}
are satisfied  with $\epsilon_n$ and $\delta_n^2$ satisfying the inequalities (\ref{ineq37}) and (\ref{ineq38}) respectively.
 Notice that  for all $ x \geq 0,$ the following three inequalities  hold
 $$\frac{x^3}{\sqrt{n}} +  \frac{x^2}{ n } +  (1+x)  \frac{ \log n }{\sqrt{n}}  \leq 2\bigg(\frac{x^3}{\sqrt{n}}  +  (1+x) \frac{ \log n }{\sqrt{n}}  \bigg),\ \ \ 0< p < 1/2, $$
  $$ \frac{x^3}{  n^{(3-4p)/2} }  + \frac{x^2}{ n } + (1+x) \frac{ \log n }{n^{(3-4p)/2} }   \leq 2\bigg(\frac{x^3}{ n^{(3-4p)/2} }  +  (1+x) \frac{ \log n }{n^{(3-4p)/2} }   \bigg),\ \ \  1/2 < p < 3/4, $$
  and
   $$ \frac{x^3}{\sqrt{ \log n}}  + \frac{x^2}{\log n }+ (1+x) \frac{   \log \log n}{\sqrt{\log n}}    \leq 2\bigg(\frac{x^3}{\sqrt{ \log n}}  +  (1+x) \frac{   \log \log n}{\sqrt{\log n}}  \bigg),\ \ \  p = 3/4. $$
Applying Lemma \ref{lemma3} to $M_n/\sqrt{v_n}$, we obtain the desired results of Theorem \ref{thmg2}.

\subsection{Proof of Corollary \ref{co03}}
 We only present the proof of Corollary \ref{co03} for $0< p <1/2$. For $1/2 < p \leq 3/4$, the proof of Corollary \ref{co03} is similar.
To prove Corollary \ref{co03}, we need the following inequalities for the normal distribution function:
\begin{eqnarray}\label{fgsgjd}
\frac{1}{\sqrt{2 \pi}(1+x)} e^{-x^2/2} \leq 1-\Phi ( x ) \leq \frac{1}{\sqrt{ \pi}(1+x)} e^{-x^2/2}, \ \ \ \   x\geq 0.
\end{eqnarray}

First we show that
\begin{eqnarray}\label{ddsmsf}
 \limsup_{n\rightarrow \infty}\frac{1}{b_n^2}\ln \mathbf{P}\bigg(\frac{   a_n S_n }{b_n \sqrt{v_n} }   \in B \bigg) \leq  - \inf_{x \in \overline{B}}\frac{x^2}{2}.
\end{eqnarray}
When $B  =\emptyset,$ the last inequality holds obviously, with   $-\inf_{x \in \emptyset} \frac{x^2}{2}=-\infty.$  Thus, we may assume that $B  \neq \emptyset.$
For a given Borel set $B\subset \mathbf{R},$ denote $x_0=\inf_{x\in B} |x|.$  Notice that $\overline{B} \supset B,$  which leads to  $x_0\geq\inf_{x\in \overline{B}} |x|$
and $x_0^2/2\geq\inf_{x\in \overline{B}}  x^2/2$.
Therefore, from  Theorem \ref{thmg2}, it follows that
\begin{eqnarray*}
\mathbf{P}\bigg(\frac{   a_n S_n }{b_n \sqrt{v_n} }   \in B \bigg)
 &\leq&  \mathbf{P}\bigg(\, \Big|\frac{   a_n S_n }{ \sqrt{v_n} }  \Big|  \geq   b_n x_0\bigg)\\
 &\leq&  2\Big( 1-\Phi \left( b_nx_0\right)\Big)
  \exp\Bigg\{ C_{p}  \bigg(  \frac{ ( b_nx_0)^3}{\sqrt{n}}    + (1+ b_nx_0)\frac{\log n}{\sqrt{n}}  \bigg) \Bigg\}.
\end{eqnarray*}
Using   (\ref{fgsgjd}),
we get
\begin{eqnarray*}
\limsup_{n\rightarrow \infty}\frac{1}{b_n^2}\ln \mathbf{P}\bigg(\frac{   a_n S_n }{b_n \sqrt{v_n} }   \in B \bigg)
 \ \leq \  -\frac{x_0^2}{2} \ \leq \  - \inf_{x \in \overline{B}}\frac{x^2}{2} ,
\end{eqnarray*}
which gives (\ref{ddsmsf}).

Next we prove that
\begin{eqnarray}\label{dfgk02}
\liminf_{n\rightarrow \infty}\frac{1}{b_n^2}\ln \mathbf{P}\bigg(\frac{   a_n S_n }{b_n \sqrt{v_n} }   \in B \bigg)  \geq   - \inf_{x \in B^o}\frac{x^2}{2} .
\end{eqnarray}
 The last inequality holds obviously for $B^o =\emptyset$, with   $-\inf_{x \in \emptyset} \frac{x^2}{2}=-\infty.$ So, we assume that $B^o \neq \emptyset$.
 Notice that $B^o$ is an open set. Therefore, for any given $\varepsilon_1>0,$ there exists an $x_0 \in B^o,$ such that
\begin{eqnarray}\label{fsdf3}
 0< \frac{x_0^2}{2} \leq   \inf_{x \in B^o}\frac{x^2}{2} +\varepsilon_1.
\end{eqnarray}
Without loss of generality, we assume that $x_0>0.$ For  all small enough $\varepsilon_2 \in (0, x_0),$ it holds $(x_0-\varepsilon_2, x_0+\varepsilon_2]  \subset B^o \subset B.$ Clearly, we have
\begin{eqnarray*}
\mathbf{P}\bigg(\frac{   a_n S_n }{b_n \sqrt{v_n} }   \in B \bigg)  &\geq&   \mathbf{P}\bigg( \frac{   a_n S_n }{  \sqrt{v_n} }  \in (b_n ( x_0-\varepsilon_2), b_n( x_0+\varepsilon_2)] \bigg)\\
&=&   \mathbf{P}\bigg(  \frac{   a_n S_n }{  \sqrt{v_n} }   \geq  b_n ( x_0-\varepsilon_2)   \bigg)-\mathbf{P}\bigg( \frac{   a_n S_n }{  \sqrt{v_n} }  \geq  b_n( x_0+\varepsilon_2) \bigg).
\end{eqnarray*}
From Theorem \ref{thmg2}, it is easy to see that
 $$\lim_{n\rightarrow \infty} \frac{\mathbf{P}\Big( \frac{   a_n S_n }{  \sqrt{v_n} }  \geq  b_n( x_0+\varepsilon_2) \Big) }{\mathbf{P}\Big(  \frac{   a_n S_n }{  \sqrt{v_n} }   \geq  b_n ( x_0-\varepsilon_2)   \Big)} =0 .$$
Using (\ref{fgsgjd}), we get
\begin{eqnarray*}
\liminf_{n\rightarrow \infty}\frac{1}{b_n^2}\ln \mathbf{P}\bigg(\frac{   a_n S_n }{b_n \sqrt{v_n} }   \in B \bigg)  \geq  -  \frac{1}{2}( x_0-\varepsilon_2)^2 . \label{ffhms}
\end{eqnarray*}
Letting $\varepsilon_2\rightarrow 0,$  by (\ref{fsdf3}), we have
\begin{eqnarray*}
\liminf_{n\rightarrow \infty}\frac{1}{b_n^2}\ln \mathbf{P}\bigg(\frac{   a_n S_n }{b_n \sqrt{v_n} }   \in B \bigg)  \ \geq\ -  \frac{x_0^2}{2}  \  \geq \   -\inf_{x \in B^o}\frac{x^2}{2} -\varepsilon_1.
\end{eqnarray*}
Because $\varepsilon_1$ can be arbitrarily small, we obtain (\ref{dfgk02}).
Combining (\ref{ddsmsf}) and (\ref{dfgk02}) together, we complete  the proof of Corollary \ref{co03} for $0< p <1/2$.

\subsection{Proof of Theorem \ref{lclt}}
  We only present the proof of point [i].  Points [ii] and [iii] can be proved in a similar way.
  We first consider the case of $1\leq k =o(n^{2/3}).$ It is easy to see that
  	\begin{eqnarray*}
\mathbf{P}( S_n=k  ) &=& \mathbf{P}( k-1< S_n \leq k  )  \nonumber \\
& =& \mathbf{P}(a_n (k-1) /\sqrt{v_n}  < a_n S_n/\sqrt{v_n} \leq a_n k /\sqrt{v_n}  )  \nonumber \\
&=&\mathbf{P}( a_n S_n/\sqrt{v_n} \leq a_n k /\sqrt{v_n}  ) - \mathbf{P}(a_n S_n/\sqrt{v_n} \leq a_n (k-1) /\sqrt{v_n}  ).
\end{eqnarray*}
For simplicity of notation, denote
$$x_k= a_n k /\sqrt{v_n},\ \ \ \ 0\leq k =o(n^{2/3}). $$
By (\ref{a10}) and (\ref{a15}), we have $x_k =o(n^{1/6})$ for all $0\leq k =o(n^{2/3}).$
From Corollary \ref{co02}, it is easy to see  that for all $0\leq k =o(n^{2/3}),$
\begin{eqnarray}
  && \mathbf{P}( a_n S_n/\sqrt{v_n} \leq x_k  ) - \Phi\big(x_k \big) \nonumber \\
  &&= (1- \Phi\big(x_k \big))- \mathbf{P}( a_n S_n/\sqrt{v_n} > x_k  ) \nonumber  \\
   &&\leq \Big(1-\Phi\big(x_k \big)  \Big) \bigg(1+ C \Big(\frac{x_k^3}{\sqrt{n}}   +  (1+x_k) \frac{ \log n }{\sqrt{n}}  \Big)\bigg) \nonumber  \\
  &&\leq \Big(1-\Phi\big(x_{k-1} \big)  \Big) \bigg(1+ C \Big(\frac{x_k^3}{\sqrt{n}}   +  (1+x_k) \frac{ \log n }{\sqrt{n}}  \Big)\bigg)  \label{dffss}
\end{eqnarray}
and
\begin{eqnarray}
&&\mathbf{P}( a_n S_n/\sqrt{v_n} \leq x_{k-1}  ) - \Phi\big(x_{k-1} \big)  \nonumber \\
 &&\geq  \Big(1-\Phi\big(x_{k-1} \big)  \Big) \bigg(1 - C \Big(\frac{x_{k-1}^3}{\sqrt{n}}   +  (1+x_{k-1}) \frac{ \log n }{\sqrt{n}}  \Big)\bigg) \nonumber  \\
&&\geq  \Big(1-\Phi\big(x_{k-1} \big)  \Big) \bigg(1 - C \Big(\frac{x_{k}^3}{\sqrt{n}}   +  (1+x_{k}) \frac{ \log n }{\sqrt{n}}  \Big)\bigg). \label{dffst}
\end{eqnarray}
Using the  inequalities (\ref{dffss}) and (\ref{dffst}), we deduce that for all $1\leq k =o(n^{2/3}),$
\begin{eqnarray}
&& \mathbf{P}( S_n=k  ) \nonumber  \\
&&=\mathbf{P}( a_n S_n/\sqrt{v_n} \leq x_k  ) - \mathbf{P}(a_n S_n/\sqrt{v_n} \leq x_{k-1 }  ) \nonumber \\
 &&\leq\Phi\big(x_k \big) -  \Phi\big(x_{k-1} \big)     +2 C  \Big(1-\Phi\big(x_{k-1} \big)  \Big)  \Big(\frac{x_k^3}{\sqrt{n}}   +  (1+x_k) \frac{ \log n }{\sqrt{n}}  \Big) \nonumber \\
&&=\int_{x_{k-1}}^{x_{k}} \frac{1}{\sqrt{2 \pi  }}\exp\Big\{-  \frac{1}{2 }t^2  \Big\} dt  +2 C  \Big(1-\Phi\big(x_{k-1} \big)  \Big)\Big(\frac{x_k^3}{\sqrt{n}}   +  (1+x_k) \frac{ \log n }{\sqrt{n}}  \Big)\nonumber \\
&&\leq \frac{1}{\sqrt{2 \pi  }}\exp\Big\{-  \frac{1}{2 }x_{k-1}^2  \Big\}  |x_k-x_{k-1} | \nonumber \\
 && +2 C  \Big(1-\Phi\big(x_{k-1} \big)  \Big)  \Big(\frac{x_k^3}{\sqrt{n}}   +  (1+x_k) \frac{ \log n }{\sqrt{n}}  \Big)   \nonumber  \\
&&\leq\frac{a_n}{\sqrt{2 \pi v_n}}\exp\Big\{-  \frac{1}{2 }x_{k-1}^2  \Big\} +2 C  \Big(1-\Phi\big(x_{k-1} \big)  \Big)  \Big(\frac{x_k^3}{\sqrt{n}}   +  (1+x_k) \frac{ \log n }{\sqrt{n}}  \Big)\label{sfgsd}.
\end{eqnarray}
Applying  (\ref{fgsgjd}) to (\ref{sfgsd}), we get for all $1\leq k =o(n^{2/3}),$
\begin{eqnarray*}
 \frac{\mathbf{P}( S_n=k  )}{ \frac{a_n}{\sqrt{2 \pi v_n}}\exp\big\{-  \frac{1}{2 }x_{k-1}^2  \big\} } &\leq& 1+   \frac{2 C   }{\sqrt{ \pi}(1+x_{k-1} )}  \Big(\frac{x_k^3}{\sqrt{n}}   +  (1+x_k) \frac{ \log n }{\sqrt{n}}  \Big) \\
 &\leq& 1+    2 C    \Big(\frac{x_k^2}{\sqrt{n}}   +   \frac{ \log n }{\sqrt{n}}  \Big).
\end{eqnarray*}
Hence, we have for all $1\leq k =o(n^{2/3}),$
\begin{eqnarray*}
 \frac{\mathbf{P}( S_n=k  )}{ \frac{a_n}{\sqrt{2 \pi v_n}}\exp\big\{-  \frac{1}{2 }x_{k}^2  \big\} }
 &\leq& \exp\Big\{ \frac{1}{2 }(x_{k}^2 -x_{k-1}^2)  \Big\} \bigg(1+    2 C    \Big(\frac{x_k^2}{\sqrt{n}}   +   \frac{ \log n }{\sqrt{n}}  \Big)\bigg) \\
 &\leq&\exp\Big\{ \frac{k a_n^2}{ v_n} \Big\} \bigg(1+    2 C   \Big(\frac{x_k^2}{\sqrt{n}}   +   \frac{ \log n }{\sqrt{n}}  \Big)\bigg).
\end{eqnarray*}
By (\ref{idf2}), we have $a_n^2/v_n \asymp 1/n$ as $n\rightarrow \infty,$ and $\exp\Big\{ \frac{k a_n^2}{ v_n} \Big\}  \leq1+ C  \frac k n$ for all $1\leq k =o(n ).$
   Thus, it holds for all $1\leq k =o(n^{2/3}),$
\begin{eqnarray}
 \frac{\mathbf{P}( S_n=k  )}{ \frac{a_n}{\sqrt{2 \pi v_n}}\exp\big\{-  \frac{1}{2 }x_{k}^2  \big\} }
  &\leq& \Big(1+ C \frac{k}{n}   \Big)\bigg(1+    2 C    \Big(\frac{x_k^2}{\sqrt{n}}   +   \frac{ \log n }{\sqrt{n}}  \Big)\bigg) \nonumber \\
  &\leq&  1+    2 C'     \Big(\frac{k}{n}+  \frac{ k^2}{ n^{3/2}}   +   \frac{ \log n }{\sqrt{n}}  \Big) .\label{hksm2}
\end{eqnarray}
Similarly, we can prove that for all $1\leq k =o(n^{2/3}),$
\begin{eqnarray}\label{hksm1}
 \frac{\mathbf{P}( S_n=k  )}{ \frac{a_n}{\sqrt{2 \pi v_n}}\exp\big\{-  \frac{1}{2 }x_{k}^2  \big\} }
   &\geq&  1 -   2 C     \Big(\frac{k}{n}+\frac{ k^2}{ n^{3/2}}   +   \frac{ \log n }{\sqrt{n}}  \Big) .
\end{eqnarray}
Combining (\ref{hksm2}) and (\ref{hksm1}) together, we get the desired inequality  of point [i]  and $1\leq k =o(n^{2/3}).$
For the case of $0\leq -k =o(n^{2/3}),$ the proof is similar.

\subsection{Proof of Corollary \ref{th025} }
We first prove the point [i] of  Corollary \ref{th025}.
By (\ref{a10}) and (\ref{a15}), we have $\frac{a_n^2}{  v_n } \asymp \frac{1}{ n} .$
Thus, from Theorem \ref{lclt}, it is easy to see that
\begin{eqnarray}
&& \sup_{|k| \leq  n^{5/8}  }  \bigg|\mathbf{P}( S_n=k    )  -    \frac{a_n}{\sqrt{2 \pi v_n}}\exp\Big\{-  \frac{(a_n k)^2}{2 v_n}  \Big\}  \bigg| \nonumber \\
 &&\leq \sup_{|k| \leq  n^{5/8}  }  C_p \frac{a_n}{\sqrt{2 \pi v_n}}\exp\Big\{-  \frac{(a_n k)^2}{2 v_n}  \Big\}  \bigg(\frac{|k|}{n}+ \frac{k^2}{ n^{3/2}} +\frac{\log n}{\sqrt{n}} \bigg) \nonumber\\
&&\leq \sup_{|k| \leq  n^{5/8}  }  C_p' \frac{1}{\sqrt{n}}\exp\Big\{-  \frac{ k^2}{n C_p'}  \Big\}  \bigg(\frac{|k|}{n}+ \frac{k^2}{ n^{3/2}} +\frac{\log n}{\sqrt{n}} \bigg) \nonumber\\
&& =  C_p' \frac{1}{n}\sup_{|k| \leq  n^{5/8}  }\exp\Big\{-  \frac{1}{ C_p'} \frac{k^2 }{n} \Big\}  \bigg( \frac{|k|}{\sqrt{n}}  + \frac{k^2}{ n} +\log n   \bigg) \nonumber\\
&& \leq   C_p''\, \frac{\log n}{n}  . \label{ineq010}
\end{eqnarray}
 Notice that
 \begin{eqnarray}
&& \sup_{|k| \geq n^{5/8}   }  \bigg|\mathbf{P}( S_n=k    )  -    \frac{a_n}{\sqrt{2 \pi v_n}}\exp\Big\{-  \frac{(a_n k)^2}{2 v_n}  \Big\}  \bigg| \nonumber \\
 &&\leq  \sup_{|k| \geq n^{5/8}   } \mathbf{P}( |S_n| \geq |k|    ) +  \sup_{|k| \geq n^{5/8}  }  \frac{a_n}{\sqrt{2 \pi v_n}}\exp\Big\{-  \frac{(a_n k)^2}{2 v_n}  \Big\} \nonumber  \\
  &&\leq    \mathbf{P}( |S_n| \geq   n^{5/8}   ) +  \sup_{|k| \geq n^{5/8}    } C_p  \frac{1 }{\sqrt{n}} \exp\Big\{-  \frac{  k ^2}{  n C_p}  \Big\}   \nonumber \\ &&\leq    \mathbf{P}( |S_n| \geq   n^{5/8}   ) +    C_p'  \frac{1 }{ n }    \nonumber .
\end{eqnarray}
Using Theorem \ref{thmg2}, we can deduce that $\mathbf{P}( |S_n| \geq   n^{5/8}   )  \leq C_p \frac{1 }{ n }. $
Thus, it holds
  \begin{eqnarray}
 \sup_{|k| \geq n^{5/8}   }  \Big|\mathbf{P}( S_n=k    )  -    \frac{a_n}{\sqrt{2 \pi v_n}}\exp\Big\{-  \frac{(a_n k)^2}{2 v_n}  \Big\}  \Big|
 \leq     C_p   \frac{1 }{ n }    \label{fds213}.
\end{eqnarray}
Combining (\ref{ineq010}) and (\ref{fds213}) together,   we obtain the desired inequality of point [i].
Point [ii]   of  Corollary \ref{th025} can be proved by a similar argument, with $n^{(9-4p)/12}$ replacing $n^{5/8}$.

\section{Conclusion}\label{conclusion}

The ERW was introduced  in order to study the memory effects in the non-Markovian random walk.
In this paper, we  study the normal approxiamtions  for the ERW, including
 Berry-Esseen's bounds, Cram\'{e}r moderate deviations and the local limit theorems.
These results can be regarded as refinements of the CLT.
Berry-Esseen's bounds (cf. Theorem  \ref{thm1}) make us better understand how the memory parameter $p$ effects the convergence rates for the absolute error of normal approximations.
From Cram\'{e}r moderate deviations, we find the domains of attraction of normal distribution for various   $p \in (0, 3/4].$ See Corollary \ref{co02}.
The local limit theorem (cf.\ Theorem  \ref{lclt}) gives  the asymptotic probability for the position of elephant.
As applications, our results can be applied to statistical inferences of the memory parameter $p$ and interval estimations for the position of the elephant.

\section*{Acknowledgements}
The authors deeply indebted to   the editor and the anonymous referees  for their helpful comments.
Fan would like to thank Quansheng Liu for his helpful suggestions.
This work has been partially supported by the National Natural Science Foundation
of China (Grant Nos.\,11601375 and 11971063).

\end{document}